\documentclass[12pt,oneside]{amsart}
\newtheorem{thm}{Theorem}[section] 
\newtheorem{cor}[thm]{Corollary}
\newtheorem{lem}[thm]{Lemma} 
\newtheorem{prop}[thm]{Proposition}
\theoremstyle{definition}
\theoremstyle{remark}
\theoremstyle{proof}

\numberwithin{equation}{section}

\newcommand{\norm}[1]{\left\Vert#1\right\Vert}
\newcommand{\abs}[1]{\left\vert#1\right\vert}
\newcommand{\set}[1]{\left\{#1\right\}}

\newcommand{\scalar}[1]{\left \langle #1 \right \rangle}
\newcommand{\sscalar}[1]{\langle #1 \rangle}
\newcommand{\Real}{\mathbb{R}}

\newcommand{\eps}{\varepsilon}
\newcommand{\To}{\longrightarrow}

\newcommand{\BP}{\mathcal{BP}}

\newcommand{\E}{\mathcal{E}}
\newcommand{\M}{\mathcal{M}}

\newcommand{\Vol}[1]{\textrm{Vol} (#1)}

\renewcommand{\S}{\mathcal{S}}

\newcommand{\RC}{\mathcal{RC}}

\begin{document}

\title{A comment on the low-dimensional Busemann-Petty problem}
\author{Emanuel Milman\\ \bigskip \fontfamily{cmr} \fontseries{m} \fontshape{sc} \fontsize{10}{0} \selectfont
Department of Mathematics, \\
The Weizmann Institute of Science, \\Rehovot 76100, Israel. \\
\medskip E-Mail: emanuel\_milman@hotmail.com.}
\thanks{Supported in part by BSF and ISF}
\begin{abstract}
The generalized Busemann-Petty problem asks whether
centrally-symmetric convex bodies having larger volume of all
$m$-dimensional sections necessarily have larger volume. When $m
> 3$ this is known to be false, but the cases $m=2,3$ are still
open. In those cases, it is shown that when the smaller body's
radial function is a $n-m$-th root of the radial function of a
convex body, the answer to the generalized Busemann-Petty problem
is positive (for any larger star-body). Several immediate
corollaries of this observation are also discussed.
\end{abstract}


\maketitle


\section{Introduction}

Let $\Vol{L}$ denote the Lebesgue measure of a set $L \subset
\Real^n$ in its affine hull, and let $G(n,k)$ denote the Grassmann
manifold of $k$ dimensional subspaces of $\Real^n$. Let $D_n$
denote the Euclidean unit ball, and $S^{n-1}$ the Euclidean
sphere. All of the bodies considered in this note will be assumed
to be centrally symmetric star-bodies, defined by a continuous
radial function $\rho_K(\theta) = \max \{r \geq 0 \; | \; r\theta
\in K \}$ for $\theta \in S^{n-1}$ and a star-body $K$.

The Busemann-Petty problem, first posed in \cite{Busemann-Petty},
asks whether two centrally-symmetric convex bodies $K$ and $L$ in
$\Real^n$ satisfying:
\begin{equation}\label{eq:Busemann-Petty}
\Vol{K\cap H} \leq \Vol{L \cap H} \; \; \forall H \in G(n,n-1)
\end{equation}
necessarily satisfy $\Vol{K} \leq \Vol{L}$. For a long time this
was believed to be true (this is certainly true for $n=2$), until
a first counterexample was given in \cite{Larman-Rogers} for a
large value of $n$. In the same year, the notion of an
\emph{intersection-body} was first introduced by Lutwak in
\cite{Lutwak-dual-mixed-volumes} (see also
\cite{Lutwak-intersection-bodies} and Section \ref{sec:two} for
definitions) in connection to the Busemann-Petty problem. It was
shown in \cite{Lutwak-intersection-bodies} (and refined in
\cite{Gardner-BP-5dim}) that the answer to the Busemann-Petty
problem is equivalent to whether all convex bodies in $\Real^n$
are intersection bodies. Subsequently, it was shown in a series of
results (\cite{Larman-Rogers}, \cite{Ball-BP}, \cite{Bourgain-BP},
\cite{Giannopoulos-BP}, \cite{Papadimitrakis-BP},
\cite{Gardner-BP-5dim}, \cite{Gardner-BP-3dim},
\cite{Koldobsky-lp-intersection-bodies},
\cite{Zhang-Correction-4dim}, \cite{GKS}), that this is true for
$n\leq 4$, but false for $n \geq 5$.

In \cite{Zhang-Gen-BP}, Zhang considered a natural generalization
of the Busemann-Petty problem, which asks whether two
centrally-symmetric convex bodies $K$ and $L$ in $\Real^n$
satisfying:
\begin{equation}\label{eq:gen-BP}
\Vol{K\cap H} \leq \Vol{L \cap H} \; \; \forall H \in G(n,n-k)
\end{equation}
necessarily satisfy $\Vol{K} \leq \Vol{L}$, where $k$ is some
integer between $1$ and $n-1$. Zhang showed that the
\emph{generalized $k$-codimensional Busemann-Petty problem} is
also naturally associated to another class of bodies, which will
be referred to as \emph{$k$-Busemann-Petty bodies} (note that
these bodies are referred to as \emph{$n-k$-intersection bodies}
in \cite{Zhang-Gen-BP} and \emph{generalized $k$-intersection
bodies} in \cite{Koldobsky-I-equal-BP}), and that the generalized
$k$-codimensional problem is equivalent to whether all convex
bodies in $R^n$ are $k$-Busemann-Petty bodies. Analogously to the
original problem, it was shown in \cite{Zhang-Gen-BP} that if $K$
and $L$ are two centrally-symmetric star-bodies (not necessarily
convex) satisfying (\ref{eq:gen-BP}), and if $K$ is a
$k$-Busemann-Petty body, then $\Vol{K} \leq \Vol{L}$.

It was shown in \cite{Bourgain-Zhang}, and later in
\cite{Koldobsky-I-equal-BP}, that the answer to the generalized
$k$-codimensional problem is negative for $k < n-3$, but the cases
$k=n-3$ and $k=n-2$ still remain open (the case $k=n-1$ is
obviously true). A partial answer to the case $k=n-2$ was given in
\cite{Bourgain-Zhang}, where it was shown that when $L$ is a
Euclidean ball and $K$ is convex and sufficiently close to $L$,
the answer is positive. Our main observation in this note
concerns the cases $k=n-2,n-3$ and reads as follows:


\begin{thm} \label{thm:main}
Let $K$ denote a centrally-symmetric convex body in $\Real^n$.
For $a=2,3$, let $K_a$ be the star-body defined by $\rho_{K_a} =
\rho_K^{1/(n-a)}$. Then $K_a$ is a $(n-a)$-Busemann-Petty body,
implying a positive answer to the $(n-a)$-codimensional
Busemann-Petty problem (\ref{eq:gen-BP}) for the pair $K_a,L$ for
any star-body $L$.
\end{thm}

The case $a=1$ is also true, but follows trivially since it is
easy to see (e.g. \cite{EMilman-GenIntBodies}) that any star-body
is an $n-1$-Busemann-Petty body. The case $a=2$ follows from $a=3$
by a general result from \cite{EMilman-GenIntBodies}, stating
that if $K$ is a $k$-Busemann-Petty body and $L$ is given by
$\rho_L = \rho_K^{k/l}$ for $1 \leq k < l \leq n-1$, then $L$ is
a $l$-Busemann-Petty body.

\medskip

Theorem \ref{thm:main} has several interesting consequences. The
first one is the following complementary result to the one
aforementioned from \cite{Bourgain-Zhang}. Roughly speaking, it
states that any small enough perturbation $K$ of the Euclidean
ball, for which we have control over the second derivatives of
$\rho_K$, satisfies the low-dimensional generalized
Busemann-Petty problem (\ref{eq:gen-BP}) with \emph{any} star-body
$L$.

\begin{cor} \label{cor:BZ}
For any $n$, there exists a function $\gamma : [0,\infty)
\rightarrow (0,1)$, such that the following holds: let $\varphi$
denote a twice continuously differentiable function on $S^{n-1}$
such that:
\[
\max_{\theta \in S^{n-1}} \abs{\varphi(\theta)} \leq 1 ,
\max_{\theta \in S^{n-1}} \abs{\varphi_{i}(\theta)} \leq M,
\max_{\theta \in S^{n-1}} \abs{\varphi_{i,j}(\theta)} \leq M,
\]
for every $i,j=1,\ldots,n-1$, where $\varphi_i$ and
$\varphi_{i,j}$ denote the first and second partial derivatives
of $\varphi$ (w.r.t. any local coordinate system of $S^{n-1}$),
respectively. Then the star-body $K^\eps$ defined by
$\rho_{K^\eps} = 1 + \eps \varphi$ for any $\abs{\eps} <
\gamma(M)$ is a $(n-a)$-Busemann-Petty body for $a=2,3$, implying
a positive answer to the $(n-a)$-codimensional Busemann-Petty
problem (\ref{eq:gen-BP}) for $K^\eps$ and any star-body $L$.
\end{cor}

Note that the definition of $K_a$ in Theorem \ref{thm:main} is
highly non-linear with respect to $K$. Since the class of
$k$-Busemann-Petty bodies is closed under certain natural
operations (see \cite{EMilman-GenIntBodies} for the latest known
results), we can take advantage of this fact to strengthen the
result of Theorem \ref{thm:main}. For instance, it is well known
(e.g. \cite{Grinberg-Zhang}, \cite{EMilman-GenIntBodies}) that the
class of $k$-Busemann-Petty bodies is closed under taking
\emph{$k$-radial sums}. The $k$-radial sum of two star-bodies
$L_1,L_2$ is defined as the star-body $L$ satisfying $\rho_L^{k}
= \rho_{L_1}^{k} + \rho_{L_2}^{k}$. When $k=1$ this operation
will simply be referred to as radial sum. The space of
star-bodies in $\Real^n$ is endowed with the natural radial metric
$d_r$, defined as $d_r(L_1,L_2) = \max_{\theta \in S^{n-1}}
\abs{\rho_{L_1}(\theta)- \rho_{L_2}(\theta)}$. We will denote by
$\RC^n$ the closure in the radial metric of the class of all
star-bodies in $\Real^n$ which are finite radial sums of
centrally-symmetric convex bodies. It should then be clear that:

\begin{cor} \label{cor:main-strengthening}
Theorem \ref{thm:main} holds for any $K \in \RC^n$.
\end{cor}

Our last remark in this note is again an immediate consequence of
Theorem \ref{thm:main} and the following characterization of
$k$-Busemann-Petty bodies due to Grinberg and Zhang
(\cite{Grinberg-Zhang}), which generalizes the characterization of
intersection-bodies (the case $k=1$) given by Goodey and Weil
(\cite{Goodey-Weil}):

\smallskip
\noindent\textbf{Theorem (Grinberg and Zhang). }\emph{A star-body
$K$ is a $k$-Busemann-Petty body iff it is the limit of
$\set{K_i}$ in the radial metric $d_r$, where each $K_i$ is a
finite $k$-radial sums of ellipsoids $\set{\E^i_j}$:
\[
\rho^k_{K_i} = \rho^k_{\E^i_1} + \ldots + \rho^k_{\E^i_{m_i}}.
\]
}

Applying Grinberg and Zhang's Theorem to the bodies $K_a$ from
Theorem \ref{thm:main}, we immediately have:

\begin{cor}
Let $K$ denote a centrally-symmetric convex body in $\Real^n$.
Then for $a=2,3$, $K$ is the limit in the radial metric $d_r$ of
star-bodies $K_i$ having the form:
\[
\rho_{K_i} = \rho^{n-a}_{\E^i_1} + \ldots +
\rho^{n-a}_{\E^i_{m_i}},
\]
where $\set{\E^i_j}$ are ellipsoids.
\end{cor}

\section{definitions and notations} \label{sec:two}

A star body $K$ is said to be an intersection body of a star body
$L$, if $\rho_K(\theta) = \Vol{L \cap \theta^\perp}$ for every
$\theta \in S^{n-1}$. $K$ is said to be an intersection body, if
it is the limit in the radial metric $d_r$ of intersection bodies
$\{K_i\}$ of star bodies $\{L_i\}$, where $d_r(K_1,K_2) =
\sup_{\theta \in S^{n-1}} \abs{\rho_{K_1}(\theta)-
\rho_{K_2}(\theta)}$. This is equivalent (e.g.
\cite{Lutwak-intersection-bodies}, \cite{Gardner-BP-5dim}) to
$\rho_K = R^*(d\mu)$, where $\mu$ is a non-negative Borel measure
on $S^{n-1}$, $R^*$ is the dual transform (as in
(\ref{eq:duality111})) to the Spherical Radon Transform
$R:C(S^{n-1}) \rightarrow C(S^{n-1})$, which is defined for $f\in
C(S^{n-1})$ as:
\[
R(f)(\theta) = \int_{S^{n-1} \cap \theta^\perp} f(\xi)
d\sigma_{n-1}(\xi),
\]
where $\sigma_{n-1}$ the Haar probability measure on $S^{n-2}$
(and we have identified $S^{n-2}$ with $S^{n-1} \cap
\theta^\perp$).

Before defining the class of $k$-Busemann-Petty bodies we shall
need to introduce the $m$-dimensional Spherical Radon Transform,
acting on spaces of continuous functions as follows:
\begin{eqnarray}
\nonumber R_m: C(S^{n-1}) \To C(G(n,m)) \\
\nonumber R_m(f) (E) = \int_{S^{n-1}\cap E} f(\theta)
d\sigma_m(\theta) ,
\end{eqnarray}
where $\sigma_m$ is the Haar probability measure on $S^{m-1}$ (and
we have identified $S^{m-1}$ with $S^{n-1} \cap E$). The dual
transform is defined on spaces of \emph{signed} Borel measures
$\M$ by:
\begin{eqnarray}
\label{eq:duality111} & R_m^*: \M(G(n,m)) \To \M(S^{n-1}) & \\
\nonumber & \int_{S^{n-1}} f R_m^*(d\mu) = \int_{G(n,m)} R_m(f)
d\mu & \forall f \in C(S^{n-1}),
\end{eqnarray}
and for a measure $\mu$ with continuous density $g$, the
transform may be explicitly written in terms of $g$ (see
\cite{Zhang-Gen-BP}):
\begin{eqnarray}
\label{eq:dual-Radon} R_m^* g (\theta) = \int_{\theta \in E \in
G(n,m)} g(E) d\nu_{m}(E) ,
\end{eqnarray}
where $\nu_{m}$ is the Haar probability measure on $G(n-1,m-1)$.

We shall say that a body $K$ is a $k$-Busemann-Petty body if
$\rho_K^k = R_{n-k}^*(d\mu)$ as measures in $\M(S^{n-1})$, where
$\mu$ is a non-negative Borel measure on $G(n,n-k)$. We shall
denote the class of such bodies by $\BP_k^n$. Choosing $k=1$, for
which $G(n,n-1)$ is isometric to $S^{n-1} / Z_2$ by mapping $H$
to $S^{n-1}\cap H^\perp$, and noticing that $R$ is equivalent to
$R_{n-1}$ under this map, we see that $\BP_1^n$ is exactly the
class of intersection bodies.

We will also require, although indirectly, several notions
regarding Fourier transforms of homogeneous distributions. We
denote by $\S(\Real^n)$ the space of rapidly decreasing
infinitely differentiable test functions in $\Real^n$, and by
$\S'(\Real^n)$ the space of distributions over $\S(\Real^n)$. The
Fourier Transform $\hat{f}$ of a distribution $f \in
\S'(\Real^n)$ is defined by $\sscalar{\hat{f},\phi} =
\sscalar{f,\hat{\phi}}$ for every test function $\phi$, where
$\hat{\phi}(y) = \int \phi(x) \exp(-i\sscalar{x,y}) dx$. A
distribution $f$ is called homogeneous of degree $p \in \Real$ if
$\scalar{f,\phi(\cdot/t)} = \abs{t}^{n+p} \scalar{f,\phi}$ for
every $t>0$, and it is called even if the same is true for
$t=-1$. An even distribution $f$ always satisfies
$(\hat{f})^\wedge = (2\pi)^n f$. The Fourier Transform of an even
homogeneous distribution of degree $p$ is an even homogeneous
distribution of degree $-n-p$.

We will denote the space of continuous functions on the sphere by
$C(S^{n-1})$. The spaces of even continuous and infinitely smooth
functions will be denoted $C_e(S^{n-1})$ and $C^\infty(S^{n-1})$,
respectively.

For a star-body $K$ (not necessarily convex), we define its
Minkowski functional as $\norm{x}_K = \min \set{ t \geq 0 \; | \;
x / t \in K}$. When $K$ is a centrally-symmetric convex body, this
of course coincides with the natural norm associated with it.
Obviously $\rho_K(\theta) = \norm{\theta}^{-1}_K$ for $\theta \in
S^{n-1}$.

\section{Proofs of the statements}

Before we begin, we shall need to recall several known facts
about the Spherical Radon Transform $R$, and its connection to
the Fourier transform of homogeneous distributions. It is well
known (e.g. \cite[Chapter 3]{Groemer}) that $R: C_e(S^{n-1})
\rightarrow C_e(S^{n-1})$ is an injective operator, and that it
is onto a dense set in $C_e(S^{n-1})$ which contains
$C_e^\infty(S^{n-1})$. The connection with Fourier transforms of
homogeneous distributions was demonstrated by Koldobsky, who
showed (e.g. \cite{Koldobsky-extremal-sections-of-lp}) the
following:

\begin{lem}
Let $L$ denote a star-body in $\Real^n$. Then for all $\theta \in
S^{n-1}$:
\[
(\norm{\cdot}_L^{-n+1})^\wedge (\theta) = \pi (n-1) \Vol{D_{n-1}}
R(\norm{\cdot}_L^{-n+1})(\theta).
\]
\end{lem}

In particular $(\norm{\cdot}_L^{-n+1})^\wedge$ is continuous, and
of course homogeneous of degree $-1$. Hence, if we denote
$\rho_K(\theta) = \norm{\theta}^{-1}_K =
(\norm{\cdot}_L^{-n+1})^\wedge(\theta)$ for $\theta \in S^{n-1}$
and use $(\norm{\cdot}^{-1}_K)^\wedge(\theta) = (2 \pi)^n
\norm{\theta}_L^{-n+1}$, we immediately get the following
inversion formula for the Spherical Radon transform:

\begin{lem} \label{lem:inversion}
Let $K$ denote a star-body in $\Real^n$ such that $\rho_K$ is in
the range of the Spherical Radon Transform. Then for all $\theta
\in S^{n-1}$:
\[
R^{-1}(\rho_K)(\theta) = \frac{ \pi (n-1) \Vol{D_{n-1}}}{(2
\pi)^n} (\norm{\cdot}_K^{-1})^\wedge(\theta).
\]
\end{lem}

Koldobsky also discovered the following property of the Fourier
transform of a norm of a \emph{convex} body (\cite[Corollary 2]
{Koldobsky-I-equal-BP}):

\begin{lem} \label{lem:whole-point}
Let $K$ be an infinitely smooth centrally-symmetric convex body
in $\Real^n$. Then for every $E \in G(n,k)$:
\[
\int_{S^{n-1} \cap E} (\norm{\cdot}_K^{-n+k+2})^\wedge(\theta)
d\theta \geq 0.
\]
\end{lem}

Since $C_e^\infty(S^{n-1})$ is in the range of the Spherical
Radon Transform, applying Lemma \ref{lem:whole-point} with
$k=n-3$ and using Lemma \ref{lem:inversion}, we have:

\begin{prop} \label{prop:whole-point}
Let $K$ be an infinitely smooth centrally-symmetric convex body
in $\Real^n$. Then for every $E \in G(n,n-3)$:
\[
\int_{S^{n-1} \cap E} R^{-1}(\rho_K)(\theta) d\theta \geq 0.
\]
\end{prop}

We are now ready to prove Theorem \ref{thm:main}.

\begin{proof}[Proof of Theorem \ref{thm:main}]
First, assume that $K$ is infinitely smooth and fix $\theta \in
S^{n-1}$. Denote by $H_\theta \in G(n,n-1)$ the hyperplane
$\theta^\perp$, and let $\sigma_{H_\theta}$ denote the Haar
probability measure on $S^{n-1} \cap H_\theta$. Let
$\eta_{H_\theta}$ denote the Haar probability measure on the
homogeneous space $G^{H_\theta}(n,n-3) := \set{E \in G(n,n-3) | E
\subset H_\theta}$, and let $\sigma_E$ denote the Haar probability
measure on $S^{n-1} \cap E$ for $E \in G(n,n-3)$. Then:
\begin{eqnarray}
\nonumber \rho_K(\theta) = R( R^{-1}(\rho_K) )(\theta) =
\int_{S^{n-1}\cap H_\theta} R^{-1}(\rho_K)(\xi) d\sigma_{H_\theta}(\xi) \\
\label{eq:plug} = \int_{E \in G^{H_\theta}(n,n-3)} \int_{S^{n-1}
\cap E} R^{-1}(\rho_K)(\xi) d\sigma_E(\xi) d\eta_{H_\theta}(E).
\end{eqnarray}
The last transition is explained by the fact that the measure
$d\sigma_E(\xi) d\eta_{H_\theta}(E)$ is invariant under orthogonal
transformations preserving $H_\theta$, so by the uniqueness of
the Haar probability measure, it must coincide with
$d\sigma_{H_\theta}(\xi)$. Denoting:
\[
g(F) = \int_{S^{n-1} \cap F^\perp} R^{-1}(\rho_K)(\xi)
d\sigma_E(\xi)
\]
for $F \in G(n,3)$, we see by Proposition \ref{prop:whole-point}
that $g \geq 0$. Plugging the definition of $g$ in
(\ref{eq:plug}), we have:
\[
\rho_K(\theta) = \int_{E \in G^{H_\theta}(n,n-3)} g(E^\perp)
d\eta_{H_\theta}(E) = \int_{F \in G_\theta(n,3)} g(F)
d\nu_{\theta}(F),
\]
where $\nu_{\theta}$ is the Haar probability measure on the
homogeneous space $G_\theta(n,3) := \set{F \in G(n,3) | \theta \in
F}$ and the transition is justified as above. By
(\ref{eq:dual-Radon}), we conclude that $\rho_K = R_3^*(g)$ with
$g \geq 0$, implying that the body $K_3$ satisfying
$\rho_{K_3}^{n-3} = \rho_K$ is in $\BP_{n-3}^n$.

As mentioned in the Introduction, the case $a=2$ follows from
$a=3$ by a general result from \cite{EMilman-GenIntBodies}, but
for completeness we reproduce the easy argument. Using
double-integration as before:
\[
\rho_K(\theta) = \int_{F \in G_\theta(n,3)} g(F) d\nu_{\theta}(F)
= \int_{J \in G_\theta(n,2)} \int_{F \in G_J(n,3)} g(F) d\nu_J(F)
d\mu_\theta(J),
\]
where $\mu_\theta$ and $\nu_J$ are the Haar probability measures
on the homogeneous spaces $G_\theta(n,2) := \set{J \in G(n,2) |
\theta \in J}$ and $G_J(n,3) := \set{F \in G(n,3) | J \subset
F}$, respectively. Denoting:
\[
h(J) = \int_{F \in G_J(n,3)} g(F) d\nu_J(F),
\]
we see that $h \geq 0$ and $\rho_K = R_2^*(h)$, implying that the
body $K_2$ satisfying $\rho_{K_2}^{n-2} = \rho_K$ is in
$\BP_{n-2}^n$.

When $K$ is a general convex body, the result follows by
approximation. It is well known (e.g. \cite[Theorem
3.3.1]{Schneider-Book}) that any centrally-symmetric convex body
$K$ may be approximated (for instance in the radial metric) by a
series of infinitely smooth centrally-symmetric convex bodies
$\set{K^i}$. Denoting by $K^i_a$ the star-bodies satisfying
$\rho_{K^i_a} = \rho_{K_i}^{1/(n-a)}$ for $a = 2,3$, we have seen
that $K^i_a \in \BP_{n-a}^n$. Obviously the series $\set{K^i_a}$
tends to $K_a$ in the radial metric, and since $\BP_{n-a}^n$ is
closed under taking radial limit (see
\cite{EMilman-GenIntBodies}), the result follows.
\end{proof}

We now turn to close a few loose ends in the proof of Corollary
\ref{cor:main-strengthening}. Since $\BP_k^n$ is closed under
$k$-radial sums, it is immediate that if $K^1$ and $K^2$ are two
convex bodies, $L$ is their radial sum, and $\rho_{T_a} =
\rho_T^{1/(n-a)}$ for $T = K_1,K_2,L$, then:
\[
\rho_{L_a}^{n-a} = \rho_L = \rho_{K_1} + \rho_{K_2} =
\rho_{K^1_a}^{n-a} + \rho_{K^2_a}^{n-a},
\]
and therefore $L_a \in \BP_{n-a}^n$. This argument of course
extends to any finite radial sum of convex bodies, and since
$\BP_k^n$ is closed under taking limit in the radial metric, the
argument extends to the entire class $\RC^n$ defined in the
Introduction.

\medskip
It remains to prove Corollary \ref{cor:BZ}.

\begin{proof}[Proof of Corollary \ref{cor:BZ}]

By Theorem \ref{thm:main}, it is enough to show that for a small
enough $\abs{\eps}$ (which depends on $n$ and $M$), the
star-bodies $L^\eps_a$ defined by $\rho_{L^\eps_a} =
\rho_{K^\eps}^{n-a}$ are in fact convex. Since $\rho_{L^\eps_a} =
(1 + \eps \varphi)^{n-a}$, it is clear that for every $\theta \in
S^{n-1}$:
\[
\abs{\rho_{L^\eps_a}(\theta)} \leq f_0(\eps,n) ,
\abs{(\rho_{L^\eps_a})_{i}(\theta)} \leq f_1(\eps,n,M),
\abs{(\rho_{L^\eps_a})_{i,j}(\theta)} \leq f_2(\eps,n,M),
\]
for every $i,j=1,\ldots,n-1$, where $f_0$ tends to 1 and
$f_1,f_2$ tend to 0, as $\eps \rightarrow 0$. It should be
intuitively clear that the convexity of $L^\eps_a$ depends only
on the behaviour of the derivatives of order 0,1 and 2 of
$\rho_{L^\eps_a}$, and since we have uniform convergence of these
derivatives to those of the Euclidean ball as $\eps$ tends to 0,
$L^\eps_a$ is convex for small enough $\eps$. To make this
argument formal, we follow \cite{Gardner-BP-problem}, and use a
formula for the Gaussian curvature of a star-body $L$ whole
radial function $\rho_L$ is twice continuously differentiable,
which was explicitly calculated in \cite[2.5]{Oliker}. In
particular, it follows that $M_L(\theta)$, the Gaussian curvature
of $\partial L$ (the hypersurface given by the boundary of $L$)
at $\rho_L(\theta) \theta$, is a continuous function of the
derivatives of order 0,1 and 2 of $\rho_L$ at the point $\theta$.
Since the Gaussian curvature of the boundary of the Euclidean
ball is a constant $1$, it follows that for small enough $\eps$,
the boundary of $L^\eps_a$ has everywhere positive Gaussian
curvature. By a standard result in differential geometry (e.g.
\cite[p. 41]{differential-geometry-book}), this implies that
$L^\eps_a$ is convex. This concludes the proof.

\end{proof}

\bibliographystyle{amsalpha}
\bibliography{../../ConvexBib}

\end{document}